\documentclass[11pt]{article}
\usepackage{amssymb}
\usepackage{amsmath}
\usepackage{amstext}
\usepackage{fancyhdr}
\usepackage{amsthm}
 \usepackage[all]{xy}

\input amssym.def
\input amssym.tex
\def\slfrac#1#2{\hbox{\kern.1em %
 \raise.5ex\hbox{\the\scriptfont0 #1}\kern-.11em %
 /\kern-.15em\lower.25ex\hbox{\the\scriptfont0 #2}}}

\newcommand{\hsp}{\hspace*{\parindent}}

\newcommand{\eeq}{\end{equation}}
\newcommand{\beql}[1]{\begin{equation}\label{#1}}

\makeatletter
\def\@sect#1#2#3#4#5#6[#7]#8{\ifnum #2>\c@secnumdepth
     \def\@svsec{}\else
     \refstepcounter{#1}\edef\@svsec{\csname the#1\endcsname.\hskip .75em }\fi
     \@tempskipa #5\relax
      \ifdim \@tempskipa>\z@
        \begingroup #6\relax
          \@hangfrom{\hskip #3\relax\@svsec}{\interlinepenalty \@M #8\par}%
        \endgroup
       \csname #1mark\endcsname{#7}\addcontentsline
         {toc}{#1}{\ifnum #2>\c@secnumdepth \else
                      \protect\numberline{\csname the#1\endcsname}\fi
                    #7}\else
        \def\@svsechd{#6\hskip #3\@svsec #8\csname #1mark\endcsname
                      {#7}\addcontentsline
                           {toc}{#1}{\ifnum #2>\c@secnumdepth \else
                             \protect\numberline{\csname the#1\endcsname}\fi
                       #7}}\fi
     \@xsect{#5}}
\def\@begintheorem#1#2{\it \trivlist \item[\hskip \labelsep{\bf #1\ #2.}]}

\def\plain{plain}\ifx\fmtname\plain\csname fi\endcsname
     
     \input docstrip
     \preamble

     Do not distribute the stripped version of this file.
     The checksum in the header refers to the documented version.

     \endpreamble
     \generateFile{here.sty}{t}{\from{here.doc}{}}
     \endinput
\fi
\ifcat a\noexpand @\let\next\relax\else\def\next{%
    \documentstyle[here,doc]{article}\MakePercentIgnore}\fi\next
\ifx\@Hxfloat\@Hundef\else\expandafter\endinput\fi
\let\@Hxfloat\@xfloat
\def\@xfloat#1[{\@ifnextchar{H}{\@HHfloat{#1}[}{\@Hxfloat{#1}[}}
\def\@HHfloat#1[H]{%
\expandafter\let\csname end#1\endcsname\end@Hfloat
\vskip\intextsep\vbox\bgroup\def\@captype{#1}\parindent\z@
\ignorespaces}

\def\end@Hfloat{\egroup\vskip \intextsep}

\makeatother

\catcode`\@=11

\renewcommand{\section}{%
        \setcounter{equation}{0}%
        \@startsection {section}{1}{\z@}{-3.5ex plus -1ex minus
        -.2ex}{2.3ex plus .2ex}{\large\bf}%
        }

\catcode`@=12

\begin{document}
\setlength{\baselineskip}{1.0\baselineskip}
\begin{center}
{\Large {\bf  Representations and cohomolgies of Rota-Baxter 3-Lie
algebras
 }}
 \\
 \vspace{1.5\baselineskip} {\large{\bf Qinxiu Sun}}
\vspace{.2\baselineskip}

Department of Mathematics, Zhejiang University of Science and Technology, Hangzhou, China, 310023\\
E-mail: qxsun@126.com

\vspace{1.5\baselineskip} {\large{\bf Shan Chen }}
\vspace{.2\baselineskip}

Department of Mathematics, Zhejiang University of Science and Technology, Hangzhou, China, 310023\\
E-mail: 1844785776@qq.com\\

\end{center}

\setlength{\baselineskip}{1.25\baselineskip}

\begin{abstract}

The goal of the present paper is to investigate representations and
 cohomologies of Rota-Baxter 3-Lie algebras with any
 weight. We introduce representations, matched pairs and Manin triples of Rota-Baxter 3-Lie algebras.
 Furthermore, we
discuss cohomology theory of Rota-Baxter 3-Lie algebras. The
deformations and central extensions of Rota-Baxter 3-Lie algebras
are also studied.
\end{abstract}

{\bf MR Subject Classification 2010}: 17B38,17B62, 16T10, 81R60.

\footnote{Keywords: 3-Lie algebra, Rota-Baxter operator,
representation, cohomology, matched pair, Manin
triple,$\mathcal{O}$-operator}

\section{Introduction}
\hsp

Rota-Baxter operators are an algebraic abstraction of the integral
operator, which originated from Baxter's study of the fluctuation
theory in probability \cite{5}. Since then, they have been related
to many topics in mathematics and mathematical physics. It has close
connection to combinatorics \cite{20, 22}, renormalization in
quantum field theory \cite{9}, multiple zeta values in number theory
\cite{24}, Yang Baxter equations \cite{1, 3}, algebraic operad
\cite{6} and so on.

Cohomologies and deformations are vital roles both in mathematics
and physics. Algebraic deformation theory arose initially in
Gerstenhaber's works for rings and algebras in \cite{18,19}.
Furthermore, Nijenhuis and Richardson considered the deformation
theory of Lie algebras \cite{28,29}. And followed by these extension
to n-Lie algebras\cite{2,10,17,30,33}. More generally, Balavoine
\cite{4} investigated deformation theory of quadratic operads.

Due to the importance of Rota-Baxter algebras, it is inevitable to
investigate their cohomologies with any weight. Recently, a great
deal of works have been devoted to this subject. The deformation and
cohomology theories of $\mathcal{O}$-operators (also named relative
Rota-Baxter operators) on Lie algebras were investigated in
\cite{34}. The authors explored deformation and cohomology theories
of relative Rota-Baxter Lie algebras of weight zero and found
applications to triangular Lie bialgebras in \cite{27}. In the same
direction, Das gave cohomology and deformation theories of
associative Rota-Baxter operators of weight 0 in \cite{11, 16} and
relative Rota-Baxter operators with any weight \cite{14}. The
cohomology and deformation of twisted Rota-Baxter operators were
considered in \cite{12,13}. Another breakthrough gained by Wang and
Zhou in \cite{36} is cohomology and homotopy theories of Rota-Baxter
algebras with any weight. Das constructed cohomologies of
Rota-Baxter operators of any weight on Lie algebras \cite{15}. Guo,
Lang and Sheng in \cite{21} introduced the notions of Rota-Baxter
Lie groups and Rota-Baxter Lie bialgebras \cite{26}, and
cohomologies of Rota-Baxter operators of weight 1 on Lie groups was
investigated in \cite{25}. There are some other works on
differential algebras with any weight \cite{23, 31, 32}.

Inspired by the previous works, we would like to study the
representation and cohomological theory of Rota-Baxter 3-Lie
algebras.

The paper is organized as follows. In Section 2, we consider
representations, matched pairs and Manin triples of Rota-Baxter
3-Lie algebras. The relation between matched pairs and Manin triples
are also characterized. In Section 3, we investigate the
cohomologies of a Rota-Baxter 3-Lie algebra with coefficients in a
representation. In Section 4, we study the deformations of
Rota-Baxter 3-Lie algebras. We also introduce the notions of a
Nijenhuis operator and $\mathcal{O}$-operator on Rota-Baxter 3-Lie
algebras. Finally, we consider the the
 central extensions of Rota-Baxter 3-Lie algebras.

\section{Representations and matched pairs of Rota-Baxter 3-Lie algebras}
A 3-Lie algebra \cite{80} is a vector space $A$ together with a
skew-symmetric linear map $[ \cdot, \cdot, \cdot]:
 \wedge^{3}A\longrightarrow A$
satisfying
$$[x_1,x_2,[y_1,y_2,y_3]]=[[x_1,x_2,y_1],y_2,y_3]+[y_1,[x_1,x_2,y_2],y_3]+[y_1,y_2,[x_1,x_2,y_3]]$$
for all $x_i,y_i\in A$.

Let $A$ be a 3-Lie algebra over field $k$
 and $\lambda\in k$. A linear operator $R:A\longrightarrow A$ is
 said to be a Rota-Baxter operator of weight $\lambda$ if it
 satisfies
\begin{eqnarray*}[R(x),R(y),R(z)]&=&R([R(x),R(y),z]+[R(x),y,R(z)]+[x,R(y),R(z)]\\&&+\lambda([R(x),y,z]+[x,R(y),z]+[x,y,R(z)])
+\lambda^{2}[x,y,z]),~~~~~~~~~~(2.1)
\end{eqnarray*} for any $x,y,z\in A$. We call $(A,
R)$ a Rota-Baxter 3-Lie algebra of weight $\lambda$.

In view of \cite{8}, we know that if $(A, R)$ is a Rota-Baxter 3-Lie
algebra of weight $\lambda$, then $(A, [ \ , \ , \ ]_{R}, R)$ is
also a Rota-Baxter 3-Lie algebra of weight $\lambda$, where $[ \ , \
, \ ]_{R}$ is given by
\begin{eqnarray*}[x,y,z]_{R}&=&[R(x),R(y),z]+[R(x),y,R(z)]+[x,R(y),R(z)]\\&&+\lambda([R(x),y,z]+[x,R(y),z]+[x,y,R(z)])
+\lambda^{2}[x,y,z],~~~~~~~~~~~~~~~~~~~~~~(2.2)
\end{eqnarray*} for any $x,y,z\in A$. Denote it by $(A_{R},R)$.

{\bf Example 2.1.} Let $({\mathfrak g}, [ \ , \ , \ ])$ be
3-dimensional 3-Lie algebra with basis ${\varepsilon_1,
\varepsilon_2, \varepsilon_3}$ and the nonzero multiplication is
given by
$$[ \varepsilon_1 , \varepsilon_2 , \varepsilon_3 ]=\varepsilon_1.$$
Define linear map $R:{\mathfrak g}\longrightarrow {\mathfrak g}$ by
the matrix $\begin {bmatrix}
 r_{11}&r_{12}&r_{13} \\
r_{21}&r_{22}&r_{23} \\
r_{31}&r_{32}&r_{33}
\end {bmatrix}$
with respect to the basis ${\varepsilon_1, \varepsilon_2,
\varepsilon_3}$. Then the linear map $R:A\longrightarrow A$ given by
the matrix
$\begin {bmatrix}
 1&a& b\\
0&1&c\\
0&d&1
\end {bmatrix}$
with respect to the basis ${\varepsilon_1,
\varepsilon_2,\varepsilon_3}$ is a Rota-Baxter operator of weight
$-1$ and the linear map $R:A\longrightarrow A$ given by the matrix
$\begin {bmatrix}
 4&a& b\\
0&1&c\\
0&d&3
\end {bmatrix}$
with respect to the basis ${\varepsilon_1,
\varepsilon_2,\varepsilon_3}$ is a Rota-Baxter operator of weight
$-2$.

A representation of a 3-Lie algebra $A$ consisting of a vector space
$V$ together with a skew-symmetric bilinear map $\rho:A\wedge A
\longrightarrow \mathfrak{gl}(V)$ satisfying
$$\rho(x_1,x_2)\rho(x_3,x_4)-\rho(x_3,x_4)\rho(x_1,x_2)
=\rho([x_1,x_2,x_3], x_4)-\rho([x_1,x_2,x_4], x_3),$$
$$\rho([x_1,x_2,x_3], x_4)=\rho(x_1,x_2)\rho(x_3,x_4)+\rho(x_2,x_3)\rho(x_1,x_4)+\rho(x_3,x_1)\rho(x_2,x_4)
,$$ for all $x_i\in A,~1\leq i\leq 4.$ We also denote
$\rho(x_1,x_2)$ by $\rho(x_1\wedge x_2)$ in the following.

One can recall the knowledge on representation of a 3-Lie algebra in
\cite{7}.\\

{\bf Definition 2.2.} Let $(A, R)$ be a Rota-Baxter 3-Lie algebra of
weight $\lambda$. A representation on Rota-Baxter 3-Lie algebra $(A,
R)$ is a triple $(V,R_{V},\rho)$, where $(V,\rho)$ is a
representation of 3-Lie algebra $A$ and satisfies
\begin{eqnarray*}&&\rho(R(x_1), R(x_2))R_{V}\\&=&R_{V}\rho(R(x_1), R(x_2))+R_{V}(\rho(R(x_1), x_2)+\rho(x_1,
R(x_2)) +\lambda \rho(x_1, x_2 ))R_{V}\\&&+\lambda
R_{V}(\rho(R(x_1), x_2)+\rho(x_1, R(x_2)) +\lambda \rho(x_1, x_2)
),~~~~~~~~~~~~~~~~~~~~~~~~~~~~~~~~~~~~~~~~ (2.3)
\end{eqnarray*}
where $x_1, x_2\in A$.\\

Let $(V,\rho)$ be a representation of 3-Lie algebra $A$. Define
$\rho^{*}:\wedge^{2}A\longrightarrow \mathfrak{gl}(V)$ by
$$\langle \rho^{*}(x,y)v^{*},w\rangle=-\langle v^{*},\rho(x_1,x_2)w\rangle$$
for any $v^{*}\in V^{*},w\in V,x,y\in A.$ Then $(V^{*},\rho^{*})$ is
again a representation of $A$. For the case of Rota-Baxter 3-Lie
algebras, we have the following result.\\

{\bf Proposition 2.3.} Let $(V,\rho,R_{V})$ be a representation of
Rota-Baxter 3-Lie algebra $(A, R)$. Then $(V^{*},\rho^{*},-\lambda
I-R_{V}^{*})$ is also a representation of $(A, R)$, where $I$ is the
identity map on $V^{*}$. We call it the dual representation.

\begin{proof} We only need to check that $(\rho^{*},-\lambda
I-R_{V}^{*})$ holds for (2.3). In fact, for any $u^{*}\in V^{*},v\in
V$, using (2.3), we get
\begin{eqnarray*}&&\langle(\lambda
I+R_{V}^{*})\rho^{*}(R(x_1), R(x_2))u^{*}-(\lambda
I+R_{V}^{*})(\rho^{*}(R(x_1), x_2)+\rho^{*}(x_1, R(x_2))
\\&&+\lambda \rho^{*}(x_1, x_2 ))(\lambda
I+R_{V}^{*})u^{*}+\lambda (\lambda I+R_{V}^{*})(\rho^{*}(R(x_1),
x_2)+\rho^{*}(x_1, R(x_2)) +\lambda \rho^{*}(x_1, x_2)
)u^{*}\\&&-\rho^{*}(R(x_1), R(x_2))(\lambda
I+R_{V}^{*})u^{*},v\rangle\\&=& \langle u^{*},\rho(R(x_1),
R(x_2))(\lambda I+R_{V})v-(\lambda I+R_{V})(\rho(R(x_1),
x_2)+\rho(x_1, R(x_2)) +\lambda \rho(x_1, x_2 ))\\&&\times(\lambda
I+R_{V})v+\lambda (\rho(R(x_1), x_2)+\rho(x_1, R(x_2)) +\lambda
\rho(x_1, x_2))(\lambda I+R_{V})v-(\lambda I+R_{V})\rho(R(X))v
\rangle\\&=&0,
\end{eqnarray*}
which yields that \begin{eqnarray*}&&(\lambda
I+R_{V}^{*})\rho^{*}(R(x_1), R(x_2))-(\lambda
I+R_{V}^{*})(\rho^{*}(R(x_1), x_2)+\rho^{*}(x_1, R(x_2)) +\lambda
\rho^{*}(x_1, x_2 ))(\lambda I+R_{V}^{*})\\&&+\lambda (\lambda
I+R_{V}^{*})(\rho^{*}(R(x_1), x_2)+\rho^{*}(x_1, R(x_2)) +\lambda
\rho^{*}(x_1, x_2) )-\rho^{*}(R(x_1), R(x_2))(\lambda
I+R_{V}^{*})=0,
\end{eqnarray*}
that is, $(\rho^{*},-\lambda I-R_{V}^{*})$ holds for (2.3).
 \end{proof}

 {\bf Example 2.4.} Let $(A,R)$ be a Rota-Baxter
3-Lie algebra of weight $\lambda$ and define ${ad}:A\wedge
A\longrightarrow A$
 by ${ad}(x_1,x_2)(x_3)=[x_1,x_2,x_3]$. Then
$(A,{ad}, R)$ is a representation of Rota-Baxter 3-Lie algebra
$(A,R)$. We call it the adjoint representation. Moreover,
$(A^{*},{ad}^{*}, -\lambda I-R^{*})$ is the dual representation, where $I$ is the identity map on $A^{*}$..\\

{\bf Proposition 2.5.} Let $(A, R)$ be a Rota-Baxter 3-Lie algebra
of weight $\lambda$ and $(V,R_{V},\rho)$ be representation of it.
Define $\tilde{\rho}:A\wedge A\longrightarrow \mathfrak{gl}(V)$ by
$$\tilde{\rho}(x_1,x_2)=\rho(R(x_1),
R(x_2))-R_{V}(\rho(R(x_1), x_2)+\rho(x_1,R(x_2)) +\lambda \rho(x_1,
x_2 )).\eqno(2.4)$$
 Then
$(V,R_{V},\tilde{\rho})$ is a representation of $(A_R,R)$. We denote
it by $\tilde{V}$.
\begin{proof}
By routine calculation, $(V,\tilde{\rho})$ is a representation of
$A_R$. We only need to check that (2.3) holds for $\tilde{\rho}$. In
fact, due to (2.3) and (2.4),
\begin{eqnarray*}&&\tilde{\rho}(R(x_1),
R(x_2))R_{V}\\&=&\rho(R^{2}(x_1),
R^{2}(x_2))R_{V}-R_{V}(\rho(R^{2}(x_1), R(x_2))+\rho(R(x_1),
R^{2}(x_2))+\lambda \rho(R(x_1),
R(x_2)))R_{V}\\&=&R_{V}\rho(R^{2}(x_1),
R^{2}(x_2))+R_{V}(\rho(R^{2}(x_1), R(x_2))+\rho(R(x_1),
R^{2}(x_2))+\lambda \rho(R(x_1), R(x_2)))R_{V}\\&&+\lambda
R_{V}(\rho(R^{2}(x_1), R(x_2))+\rho(R(x_1), R^{2}(x_2))+\lambda\rho(
R(x_1), R(x_2)))\\&&- R_{V}^{2}\rho(R^{2}(x_1), R(x_2))
-R_{V}^{2}(\rho(R^{2}(x_1), x_2)+\rho(R(x_1), R(x_2))+\lambda\rho(
R(x_1), x_2))R_{V}\\&&- \lambda R_{V}^{2}(\rho(R^{2}(x_1),
x_2)+\rho(R(x_1), R(x_2))+\lambda\rho( R(x_1), x_2)) -
R_{V}^{2}\rho(R(x_1), R^{2}(x_2))\\&&-
R_{V}^{2}(\rho(R(x_1),R(x_2)+\rho(x_1, R^{2}(x_2))+\lambda\rho( x_1,
R(x_2)))R_{V}- \lambda R_{V}^{2}(\rho(R(x_1), R(x_2))\\&&+\rho(x_1,
R^{2}(x_2))+\lambda \rho(x_1, R(x_2)))-\lambda R_{V}^{2}\rho(R(x_1),
R(x_2)) -\lambda R_{V}^{2}(\rho(R(x_1), x_2)\\&&+\rho(x_1,
R(x_2))+\lambda \rho(x_1, x_2))R_{V} -\lambda^{2}
R_{V}^{2}(\rho(R(x_1), x_2)+\rho(x_1, R(x_2))+\lambda \rho(x_1,
x_2))\\&=& R_{V}\rho(R^{2}(x_1),
R^{2}(x_2))-R_{V}^{2}(\rho(R^{2}(x_1), R(x_2))+\rho(R(x_1),
R^{2}(x_2))+\lambda \rho(R(x_1), R(x_2)))\\&&+R_{V}\rho(R^{2}(x_1),
R(x_2))- R_{V}^{2}(\rho(R^{2}(x_1), x_2)+\rho(R(x_1),
R(x_2))+\lambda \rho( R(x_1), x_2))R_{V}\\&&+R_{V}\rho(R(x_1),
R^{2}(x_2))- R_{V}^{2}(\rho(x_1, R^{2}(x_2))+\rho(R(x_1),
R(x_2))+\lambda \rho( x_1, R(x_2)))R_{V}\\&&+\lambda
R_{V}\rho(R(x_1), R(x_2))-\lambda^{2} R_{V}^{2}(\rho(R(x_1),
x_2)+\rho(x_1, R(x_2))+\lambda \rho(x_1, R(x_2)))R_{V}
\\&=&R_{V}\tilde{\rho}(R(x_1),
R(x_2))+R_{V}(\tilde{\rho}(R(x_1), x_2)+\tilde{\rho}(x_1,
R(x_2))+\lambda \tilde{\rho}(x_1, x_2))R_{V}\\&&+\lambda
R_{V}(\tilde{\rho}(R(x_1), x_2)+\tilde{\rho}(x_1, R(x_2))+\lambda
\tilde{\rho}( x_1, x_2)).
\end{eqnarray*}

\end{proof}

The above representation $\tilde{V}$ can be applied to cohomologies
of Rota-Baxter operators
and Rota-Baxter 3-Lie algebras in the next Section. W e can also construct another representation of $(A_R,R)$ from the
 given representation $(V,\rho)$ of Rota-Baxter 3-Lie algebra $(A,R)$. \\

 {\bf Proposition 2.6.} Let $(A, R)$ be a
Rota-Baxter 3-Lie algebra of weight $\lambda$ and $(V,R_{V},\rho)$
be representation of it. Define bilinear map $\bar{\rho}:A\wedge
A\longrightarrow \mathfrak{gl}(V)$ by
\begin{eqnarray*}\bar{\rho}(x_1,x_2)&=&\rho(R(x_1),
R(x_2))+(\rho(R(x_1), x_2)+\rho(x_1,R(x_2) )+\lambda \rho(x_1, x_2
))R_{V}\\&&+\lambda(\rho(R(x_1), x_2)+\rho(x_1,R(x_2) )+\lambda
\rho(x_1, x_2 )).\end{eqnarray*}
 Then $(V,R_{V},\bar{\rho})$ is a
representation of $(A_R,R)$. We denote it by $\bar{V}$.

\begin{proof}
The proof is similar to Proposition 2.5.
\end{proof}

 {\bf Definition 2.7.} Let $(A,R_{A})$ and $(B,R_{B})$
be two Rota-Baxter $3$-Lie algebras of weight $\lambda$ such that
$(B , \rho,R_{B})$ is a representation of $(A,R_{A})$ and $(A,
\varrho,R_{A})$ is a representation of $(B,R_{B})$ with bilinear
maps $\rho:\wedge^{2} A\longrightarrow \mathfrak{gl}(B)$ and
$\varrho:\wedge^{2} B\longrightarrow \mathfrak{gl}(A)$. Moreover,
$(A,B , \rho, \varrho)$ is a matched pair of 3-Lie algebras with the
3-Lie algebra structure on $A\bowtie B$ given by
\begin{eqnarray*}[x+a,y+b,z+c]&=&[x,y,z]+\rho(a,b)z+\rho(b,c)x+\rho(c,a)y+[a,b,c]
\\&&+\varrho(x,y)c+\varrho(y,z)a+\varrho(z,x)b.
\end{eqnarray*}
 Then $(A\bowtie
B,R_{A}+R_{B})$ is a Rota-Baxter $3$-Lie algebra with
 $$(R_{A}+R_{B})(x+a)=R_{A}(x)+R_{B}(a)$$ for any $x,y\in A,a,b\in B$. We call
$(A,B,R_{A},R_{B}, \rho,\varrho)$ the matched pair of Rota-Baxter $3$-Lie algebras.\\

{\bf Proposition 2.8.} Let $(A,R_{A})$ and $(B,R_{B})$ be two
Rota-Baxter $3$-Lie algebras of weight $\lambda$, $(B , \rho)$ is a
representation of $A$ and $(A, \varrho)$ is a representation of $B$
with bilinear maps $\rho:\wedge^{2} A\longrightarrow {gl}(B)$ and
$\varrho:\wedge^{2} B\longrightarrow {gl}(A)$. Then the following
conditions are equivalent

(i) $(A,B,R_{A},R_{B}, \rho, \varrho)$ is a matched pair of
Rota-Baxter $3$-Lie algebras.

(ii) $(A, B , \rho, \varrho)$ is a matched pair of 3-Lie algebras
and satisfying
\begin{eqnarray*}&&\rho(R_{A}(x_1),R_{A}(x_2))R_{B}\\&=&R_{B}(R_{A}(x_1),R_{A}(x_2))+R_{B}\rho(R_{A}(x_1), x_2)+\rho(x_1,
R_{A}(x_2)) +\lambda \rho(x_1, x_2 ))R_{B}\\&&+\lambda
R_{B}(\rho(R_{A}(x_1), x_2)+\rho(x_1, R_{A}(x_2)) +\lambda \rho(x_1,
x_2) ),~~~~~~~~~~~~~~~~~~~~~~~~~~~~~ (2.5)
\end{eqnarray*}
\begin{eqnarray*}&&\varrho(R_{B}(a_1),R_{B}(a_2))R_{A}\\&=&R_{A}(\varrho(R_{B}(a_1),R_{B}(a_2))+R_{A}\varrho(R_{B}(a_1), a_2)+\varrho(a_1,
R_{B}(a_2)) +\lambda \varrho(a_1, a_2 ))R_{A}\\&&+\lambda
R_{A}(\varrho(R_{B}(a_1), a_2)+\varrho(a_1, R_{B}(a_2)) +\lambda
\varrho(a_1, a_2) ),~~~~~~~~~~~~~~~~~~~~~~~~~~~~ (2.6)
\end{eqnarray*}
for any $x_1,x_2\in A,a_1,a_2\in B$.

(iii) $(A, B , \rho, \varrho)$ is a matched pair of 3-Lie algebras
and $$R_A+R_B:A\bowtie B\longrightarrow A\bowtie B$$ is a
Rota-Baxter operator of weight $\lambda$ on the $3$-Lie algebras
$A\bowtie B$ with $$(R_A+R_B)(x+a)=R_A(x)+R_B(a).$$

\begin{proof}
In view of the definitions of representations and matched pairs of
Rota-Baxter $3$-Lie algebras, we get that (i) is equivalent to
(iii). Subsequently, we will prove that $(iii)$ is equivalent to
(ii).

If $(A, B , \rho, \varrho)$ is a matched pair of 3-Lie algebras,
then $R_A+R_B$ is a Rota-Baxter operator of weight $\lambda$ on the
$3$-Lie algebra $A\bowtie B$ if and only if
\begin{eqnarray*}&&[(R_{A}+R_B)(x+a),(R_{A}+R_B)(y+b),(R_{A}+R_B)(z+c)]\\&=&(R_A+R_B)([(R_{A}+R_B)(x+a),(R_{A}+R_B)(y+b),z+c]
\\&&+[(R_{A}+R_B)(x+a),y+b,(R_{A}+R_B)(z+c)]\\&&+[x+a,(R_{A}+R_B)(y+b),(R_{A}+R_B)(z+c)]+\lambda[(R_{A}+R_B)(x+a),y+b,z+c]\\&&+\lambda[x+a,(R_{A}+R_B)(y+b),z+c]
+\lambda[x+a,y+b,(R_{A}+R_B)(z+c)]\\&&+\lambda^{2}[x+a,y+b,z+c]).~~~~~~~~~~~~~~~~~~~~~~~~~~~~~~~~~~~~~~~~~~~~~~~~~~~~~~~~~~~~~~~~~~~~~~~(2.7)
\end{eqnarray*}
By direct calculation, (2.7) holds if and only if (2.5), (2.6) hold.
It follows that $(iii)$ is equivalent to (ii).
\end{proof}

{\bf Proposition 2.9.} Let $(A,B,R_{A},R_{B}, \rho, \varrho)$ be a
matched pair of Rota-Baxter $3$-Lie algebras. Then $(A_R,B_R,
\bar{\rho}, \bar{\varrho})$ is a matched pair of $3$-Lie algebras.
Furthermore, $A_R\bowtie B_R\simeq (A\bowtie B)_{R}$ as $3$-Lie
algebras.
\begin{proof}
The $3$-Lie algebra structure on $(A\bowtie B)_{R}$ is given by
\begin{eqnarray*}&&[x+a,y+b,z+c]_{R}\\&=&[R_{A}x+R_{B}a,R_{A}y+R_{B}b,z+c]+[R_{A}x+R_{B}a,y+b,R_{A}z+R_{B}c]\\&&+[x+a,R_{A}y+R_{B}b,,R_{A}z+R_{B}c]+
\lambda([R_{A}x+R_{B}a,y+b,z+c]+[x+a,R_{A}y+R_{B}b,z+c]\\&&+[x+a,y+b,R_{A}z+R_{B}c]+\lambda[x+a,y+b,z+c])\\&=&
[R_{A}x,R_{A}y,z]+\rho(R_{B}a,R_{B}b)z+\rho(R_{B}b,c)R_{A}x+\rho(c,R_{B}a)R_{A}y+
[R_{B}a,R_{B}b,c]\\&&+\varrho(R_{A}x,R_{A}y)c+\varrho(R_{A}y,z)R_{B}a+\varrho(z,R_{A}x)R_{B}b+
[R_{A}x,y,R_{A}z]+\rho(R_{B}a,b)R_{A}z\\&&+\rho(b,R_{B}c)R_{A}x+\rho(R_{B}c,R_{B}a)y+
[R_{B}a,b,R_{B}c]+\varrho(R_{A}x,y)R_{B}c+\varrho(y,R_{A}z)R_{B}a\\&&+\varrho(R_{A}z,R_{A}x)b
+
[x,R_{A}y,R_{A}z]+\rho(a,R_{B}b)R_{A}z+\rho(R_{B}b,R_{B}c)x\\&&+\rho(R_{B}c,a)R_{A}y+
[a,R_{B}b,R_{B}c]+\varrho(x,R_{A}y)R_{B}c+\varrho(R_{A}y,R_{A}z)a+\varrho(R_{A}z,x)R_{B}b
\\&&+
\lambda([R_{A}x,y,z]+\rho(R_{B}a,b)z+\rho(b,c)R_{A}x+\rho(c,R_{B}a)y+
[R_{B}a,b,c]\\&&+\varrho(R_{A}x,y)c+\varrho(y,z)R_{B}a+\varrho(z,R_{A}x)b
+ [x,R_{A}y,z]+\rho(a,R_{B}b)z+\rho(R_{B}b,c)x\\&&+\rho(c,a)R_{A}y+
[a,R_{B}b,c]+\varrho(x,R_{A}y)c+\varrho(R_{A}y,z)a+\varrho(z,x)R_{B}b
+ [x,y,R_{A}z]\\&&+\rho(a,b)R_{A}z+\rho(b,R_{B}c)x+\rho(R_{B}c,a)y+
[a,b,R_{B}c]+\varrho(x,y)R_{B}c+\varrho(y,R_{A}z)a\\&&+\varrho(R_{A}z,x)b
+
\lambda[x,y,z]+\lambda\rho(a,b)z+\lambda\rho(b,c)x+\lambda\rho(c,a)y+
\lambda[a,b,c]\\&&+\lambda\varrho(x,y)c+\lambda\varrho(y,z)a+\lambda\varrho(z,x)b)\\&=&[x,y,z]_R+
\bar{\rho}(a,b)z+\bar{\rho}(b,c)x+\bar{\rho}(c,a)y+[a,b,c]_R+\bar{\varrho}(x,y)c+\bar{\varrho}(y,z)a+\bar{\varrho}(z,x)b,
\end{eqnarray*}
which implies that $(A_R,B_R, \bar{\rho}, \bar{\varrho})$ is a
matched pair of $3$-Lie algebras. Furthermore, $A_R\bowtie B_R\simeq
(A\bowtie B)_{R}$ as $3$-Lie algebras.
\end{proof}

{\bf Proposition 2.10.} Let $(A,B,R_{A},R_{B}, \rho, \varrho)$ be a
matched pair of Rota-Baxter $3$-Lie algebras of weight $\lambda$.
Then $(A_R,B_R,R_{A},R_{B}, \bar{\rho}, \bar{\varrho})$ is a matched
pair of Rota-Baxter $3$-Lie algebras of weight $\lambda$.
\begin{proof}
By Proposition 2.9, $A_R\bowtie B_R\simeq (A\bowtie B)_{R}$. Thus,
$R_A+R_B$ is a Rota-Baxter operator of $A_R\bowtie B_R$.  By
Proposition 2.8,  $(A_R,B_R,R_{A},R_{B}, \bar{\rho}, \bar{\varrho})$
is a matched pair of Rota-Baxter $3$-Lie algebras of weight
$\lambda$.
\end{proof}

{\bf Definition 2.11.} Let $(A,R)$ be a Rota-Baxter $3$-Lie algebra
of weight $\lambda$. A bilinear form $B$ on $(A,R)$ is called
invariant if it satisfies
$$B([x_1,x_2,x_3],x_4)+B([x_1,x_2,x_4],x_3)$$
and
$$B(R(x),y)+B(x,R(y)+\lambda B(x,y)=0.$$
 A Rota-Baxter $3$-Lie algebra is called pseudo-metric if there is
 a nondegenerate symmetric invariant bilinear form on $(A,R)$.\\

{\bf Definition 2.12.} A Manin triple of Rota-Baxter 3-Lie algebras
of weight $\lambda$ consists of a pseudo-metric Rota-Baxter $3$-Lie
algebra $(A,R)$ and Rota-Baxter 3-Lie algebras $(A_1,R_1)$,
$(A_2,R_2)$ such that

(i) $(A_1,R_1)$ and $(A_2,R_2)$ are Rota-Baxter Lie subalgebras of
$(A,R)$.

(ii) $A=A_1\oplus A_2$ as vector spaces;

(iii) For all $x_1,x_2\in A_1,~a_1,a_2\in A_2$, we have
$P_1[x_1,x_2,a_1]=0$ and $P_2[a_1,a_2,x_1]=0$ with $P_1,P_2$ are
projections from $A_1\oplus A_2$ to $A_1$ and $A_2$ respectively.\\

More details on matched pairs and Manin triples of 3-Lie algebras
can be found in \cite{7}.\\

 {\bf Proposition 2.13.} Let $(A,R)$ and
$(A^{*}, -\lambda I-R^{*})$ be Rota-Baxter 3-Lie algebras of weight
$\lambda$. Then $(A\oplus A^{*},B,A,A^{*} )$ is a Manin triple of
Rota-Baxter 3-Lie algebras if and only if
$(A,A^{*},ad^{*},\mathfrak{ad}^{*},R,-\lambda I-R^{*})$ is a matched
pair of Rota-Baxter 3-Lie algebras, where
$$B(x+\xi,y+\eta)=\langle x,\eta\rangle + \langle \xi,y\rangle. $$
\begin{proof}
According to Proposition 4.7 \cite{7}, we easily get the result.
\end{proof}

\section{Cohomology of Rota-Baxter 3-Lie algebras }
\hsp In this section, we first recall the cohomogies of 3-Lie
algebras \cite{2,30}. Then, we combine the cohomologies of
  3-Lie algebras $A$ and $A_R$ to construct
 cohomologies of Rota-Baxter 3-Lie
algebras.

Let $A$ be a $3$-Lie algebra and $L=\wedge^{2}A$ be the associated
Leibniz algebra. Suppose that $(\rho,V)$ is a representation of $A$,
the space $C_{3-{\hbox{Lie}}}^{p}(A,V)$ of $p$-cochains ($p\geq 1$)
is the set of multilinear maps of the form
$$f:\wedge^{p-1}L\wedge A\longrightarrow V,$$
and the coboundary operator
$\partial:C_{3-{\hbox{Lie}}}^{p}(A,V)\longrightarrow
C_{3-{\hbox{Lie}}}^{p+1}(A,V)$ is as follows:
\begin{eqnarray*}&&(\partial f)(X_1,\cdot\cdot\cdot,X_{p},z)\\&=&\sum_{1\leq i<
k\leq {p}}
(-1)^{i}f(X_1,\cdot\cdot\cdot,\hat{X_i},\cdot\cdot\cdot,X_{k-1},[X_i,X_k]_F,X_{k+1},\cdot\cdot\cdot,X_{p}
,z)
\\&&+\sum_{i=1}^{p}(-1)^{i}f(X_1,\cdot\cdot\cdot,\hat{X_i},\cdot\cdot\cdot,X_{p},[X_i,z])+
\sum_{i=1}^{p}(-1)^{i+1}\rho(X_i)f(X_1,\cdot\cdot\cdot,\hat{X_i},\cdot\cdot\cdot,X_{p},z)\\&&
+(-1)^{p+1}(\rho(y_{p},z)f(X_1,\cdot\cdot\cdot,X_{p-1},x_{p})+\rho(z,x_{p})f(X_1,\cdot\cdot\cdot,X_{p-1},y_{p}))
\end{eqnarray*}
for all $X_i=x_{i}\wedge y_{i}\in L $ and $z\in A$.

In the following, we are ready to define cohomology of Rota-Baxter
$3$-Lie algebras.

Let $(A,R)$ be a Rota-Baxter $3$-Lie algebra of any weight $\lambda$
and $( V,\rho,R_{V})$ be its representation. It is known that
$(A_{R},[ \ , \ , \ ]_R,R)$ is also a $\lambda$-weighted Rota-Baxter
$3$-Lie algebra and $(V,\tilde{\rho},R_{V})$ is a representation of
$(A_R,R)$ by Proposition 2.5.

It is natural to study the cohomologies of $3$-Lie algebra $(A_R,[ \
, \ , \ ]_R)$ with coefficients in $(\tilde{V},\tilde{\rho})$.

The coboundary operator
$\partial_{R}:C_{3-{\hbox{Lie}}}^{p}(A_R,\tilde{V})\longrightarrow
C_{3-{\hbox{Lie}}}^{p+1}(A_R,\tilde{V})$ is given by:
\begin{eqnarray*}&&(\partial_{R} f)(X_1,\cdot\cdot\cdot,X_{p},z)\\&=&\sum_{1\leq i<
k\leq {p}}
(-1)^{i}f(X_1,\cdot\cdot\cdot,\hat{X_i},\cdot\cdot\cdot,X_{k-1},[X_i,X_k]_F,X_{k+1},\cdot\cdot\cdot,X_{p}
,z)
\\&&+\sum_{i=1}^{p}(-1)^{i}f(X_1,\cdot\cdot\cdot,\hat{X_i},\cdot\cdot\cdot,X_{p},[X_i,z]_{R})+
\sum_{i=1}^{p}(-1)^{i+1}\tilde{\rho}(X_i)f(X_1,\cdot\cdot\cdot,\hat{X_i},\cdot\cdot\cdot,X_{p},z)\\&&
+(-1)^{p+1}(\tilde{\rho}(y_{p},z)f(X_1,\cdot\cdot\cdot,X_{p-1},x_{p})+\tilde{\rho}(z,x_{p})f(X_1,\cdot\cdot\cdot,X_{p-1},y_{p}))
\end{eqnarray*}

for all $X_i=x_{i}\wedge y_{i}\in L $ and $z\in A$.

We obtain a complex $(C_{\hbox{3-Lie}}^{*}(A_R,\tilde{V}),
\partial_R)$. Denote the set of closed $n$-cochains by
$\mathcal{Z}^{n}_{\hbox{3-Lie}}(A_R, \tilde{V})$ and the set of
exact $n$-cochains by $\mathcal{B}^{n}_{\hbox{3-Lie}}(A_R,
\tilde{V})$. We define the corresponding cohomology group by
$$\mathcal{H}^{n}_{\hbox{3-Lie}}(A_R,
\tilde{V})=\mathcal{Z}^{n}_{\hbox{3-Lie}}(A_R,
\tilde{V})/\mathcal{B}^{n}_{\hbox{3-Lie}}(A_R, \tilde{V}).$$

{\bf Definition 3.1.} Let $(A,R)$ be a Rota-Baxter $3$-Lie algebra
of weight $\lambda$ and $( V,\rho,R_{V})$ be its representation.
Then the cochain complex
$(C_{\hbox{3-Lie}}^{*}(A_R,\tilde{V}),\partial_{R})$ is called the
cochain complex of Rota-Baxter operator $R$ with coefficients in
$V$, the corresponding cohomology group
$\mathcal{H}^{n}_{\hbox{3-Lie}}(A_R, \tilde{V})$ is called the $n$-
cohomology of Rota-Baxter operator $R$ with coefficients in $V$.\\

{\bf Proposition 3.2.} Define linear map
$\delta:C_{\hbox{3-Lie}}^{n+1}(A,V)\longrightarrow
C_{\hbox{3-Lie}}^{n+1}(A_R,\tilde{V})$ by
$$\delta
f(x_1,\cdot\cdot\cdot,x_{2n+1})=f(R(x_1),\cdot\cdot\cdot,R(x_{2n+1}))-R_{V}
\sum_{k=0}^{2n}\lambda^{2n-k}f^{k}(x_1, \cdot\cdot\cdot,x_{2n+1})
$$
where $$f^{k}(x_1,\cdot\cdot\cdot,x_n) =
f(\underbrace{I\otimes\cdot\cdot\cdot \otimes
R\otimes\cdot\cdot\cdot \otimes I}_{R~~ \hbox{appears} ~~k
~~\hbox{times}})(x_1,\cdot\cdot\cdot,x_n).$$ Then $\delta$ is a
cochain map from the cochain complex $C_{\hbox{3-Lie}}^{*}(A,V)$ to
$ C_{\hbox{3-Lie}}^{*}(A_R,\tilde{V})$, that is, the following
diagram is commutative:
$$\xymatrix{
 C_{\hbox{3-Lie}}^{n}(A,V) \ar[d]_{\partial} \ar[r]^{\delta} & C_{\hbox{3-Lie}}^{n}(A_R,\tilde{V}) \ar[d]^{\partial_{R}} \\
  C_{\hbox{3-Lie}}^{n+1}(A,V)\ar[r]^{\delta} & C_{\hbox{3-Lie}}^{n+1}(A_R,\tilde{V})
  .}$$

\begin{proof} See the Appendix.
\end{proof}

Define the set of $(n+1)$-cochains by
  $$C_{\hbox{RB}}^{n+1}(A,V)=C_{\hbox{3-Lie}}^{n+1}(A,V)\times C_{\hbox{3-Lie}}^{n}(A_R,\tilde{V}),~(n\geq 1)$$
  and $$C_{\hbox{RB}}^{1}(A,V)=C_{\hbox{3-Lie}}^{1}(A,V).$$

Given the linear map
$\partial_{RB}:C_{\hbox{RB}}^{n}(A,V)\longrightarrow
C_{\hbox{RB}}^{n+1}(A,V)$ by
$$\partial_{RB}(f,g)=(\partial f,\partial_R
g+(-1)^{n}\delta f) ~~\hbox{for~~ any} ~~(f,g)\in
C_{\hbox{RB}}^{n}(A,V).$$

In view of Proposition 3.2, we have \\

{\bf Theorem 3.3.} The operator $\partial_{RB}$ is a coboundary
operator, that is,
$\partial_{RB}\partial_{RB}=0$.\\

Associated to the representation $(V,R_{V},\rho)$, we obtain a
cochain complex $(C_{\hbox{RB}}^{*}(A,V),
\partial_{RB})$. Denote the cohomology group of this cochain
complex by $H^{*}_{RB}(A, V)$, which is called the cohomology of the
Rota-Baxter 3-Lie algebra $(A,R)$ with coefficients in the
representation $(V,R_{V},\rho)$.

In view of the definitions of $C_{\hbox{3-Lie}}^{n}(A,V)$ and
$C_{\hbox{3-Lie}}^{n}(A_R,\tilde{V})$, we easily get an exact
sequence of cochain complexes,
$$0\longrightarrow
C_{\hbox{3-Lie}}^{n}(A_R,\tilde{V})\longrightarrow
C_{\hbox{RB}}^{n+1}(A,V) \longrightarrow
C_{\hbox{3-Lie}}^{n+1}(A,V)\longrightarrow 0.$$

{\bf Remark 3.4.} If $A$ is a 3-Lie algebra, then $(L(A)=A\wedge A,
[ \ , \ ]_{F} )$ is the associated Leibniz algebra \cite{30} with $[
\ , \ ]_{F} $ given by $$[ X, Y ]_{F}=[x_1,x_2,y_1]\wedge
y_2+y_1\wedge[x_1,x_2,y_2]$$ for any $X=x_1\wedge x_2,Y=y_1\wedge
y_2\in L(A)$. Moreover, the cohomology complex of a 3-Lie algebra
$A$ with coefficients in $V$ coincides with the Loday-Pirashvili
cohomology complex of Leibniz algebra $L(A)$ with coefficients in
$\hbox{Hom}(A, V )$ with the same coboundary operator \cite{10}.
But, in the case of Rota-Baxter 3-Lie algebras, we can not get a
similar result, since it is difficult to gain a Rota-Baxter operator
on $L(A)=A\wedge A$.

\section{Deformation of Rota-Baxter 3-Lie algebras}

Let $(A,\pi,R
 )$ be a Rota-Baxter 3-Lie algebra,
 $\pi_{i}:\wedge^{3}A
 \longrightarrow A
 $ be trilinear map and $R_{i}:A\longrightarrow
A$ be linear map.
 Consider the space $A[[t]]$ of formal power series in $t$ with
 coefficients in $A$ and a $t$-parametrized family of trilinear operations
 $$\pi_{t}(x,y,z)=\sum_{i=0}^{n-1}t^{i}\pi_{i}(x,y,z),\eqno (4.1)$$
 and linear operations
 $$R_{t}(x)=\sum_{i=0}^{n-1}t^iR_{i}(x),\eqno (4.2)$$
 where $R_{0}=R$ and $\pi_{0}=\pi$.

 If all $(A[[t]],
 \pi_{i},R_{i})$ are Rota-Baxter 3-Lie algebras, we say that
 $(\pi_{i},R_i)~(i=0,1,\cdot\cdot\cdot,)$ generate a deformation of the Rota-Baxter 3-Lie algebra $(
 A,\pi,R)$.

If all $(A[[t]],
 \pi_{i},R_{i})$ are Rota-Baxter 3-Lie algebras, then we have
\begin{eqnarray*}&&\sum_{i+j=k,~i,~j\geq
 0}(\pi_i(\pi_j(v,w,x),y,z)+\pi_i(x,\pi_j(v,w,y),z)+\pi_i(x,y,\pi_j(v,w,z))\\&&-\pi_i(v,w,\pi_j(x,y,z)))=0, ~~~~~~~~~
 ~~~~~~~~~~~~~~~~~~~~~~~~~~~~~~~~~~~~~~~~~~~~~~~~~~~~~~~~~~~~~(4.3)
\end{eqnarray*} and
\begin{eqnarray*}&&\sum_{i+j+k+l=n,~i,~j,~k,~l\geq
 0}\pi_i(R_j(x),R_k(y),R_l(z))\\&=&\sum_{i+j+k+l=n,~i,~j,~k,~l\geq
 0}R_{i}(\pi_j(R_k(x),R_l(y),z)+\pi_j(x,R_k(y),R_l(z))+\pi_j(R_k(x),y,R_l(z)))
\\&&+\lambda\sum_{i+j+k=n,~i,~j,~k\geq
 0}R_{i}(\pi_j(R_k(x),y,z)+\pi_j(x,R_k(y),z)+\pi_j(x,y,R_k(z)))\\&&+\lambda^{2}\sum_{i+j=n,~i,~j\geq
 0}R_{i}\pi_j(x,y,z).~~~~~~~~~~~~~~~~~~~~~~~~~~~~~~~~~~~~~~~~~~~~~~~~~~~~~~~~~~~~~~~~~(4.4)\end{eqnarray*}

 {\bf Proposition 4.1.}  Let $(A[[t]],\pi_{i},R_i)$ be a
 one-parameter deformation of the Rota-Baxter 3-Lie algebra
 $(A,\pi,R
 )$. Then $(\pi_1,R_1)$ is a 2-cocycle in the cohomology of the Rota-Baxter 3-Lie algebra
 $(A,\pi,R
 )$ with coefficients in itself.
\begin{proof} If $(\pi_{i},R_i)$ is a
 one-parameter deformation of the Rota-Baxter 3-Lie algebra $(A,\pi,R)$, then $(4.3)$ and $(4.4)$ hold.
 On the other hand, if $(\pi_1,R_1)$ is a 2-cocycle, then $\partial \pi_1=0$
and $\partial_{R}R_1+\delta \pi_1=0$.

On the basis of deformation theory of 3-Lie algebras \cite{38,39},
$\partial \pi_1=0$ is equivalent to (4.3) when $k=1$.

By direct calculation,
\begin{eqnarray*}&&\partial_{R}R_{1}(x,y,z)\\&=&-R_{1}([x,y,z]_{R})+\tilde{\rho}(x,y)R_{1}(z)+\tilde{\rho}(y,z)R_{1}(x)+\tilde{\rho}(z,x)R_{1}(y)
\\&=&-R_{1}([x,y,z]_{R})+(\rho(R(x),R(y))-R(\rho(R(x),y)+\rho(x,R(y))+\lambda\rho(x,y)))R_{1}(z)\\&&
+(\rho(R(y),R(z))-R(\rho(R(y),z)+\rho(y,R(z))+\lambda\rho(y,z)))R_{1}(x)
\\&&+(\rho(R(z),R(x))-R(\rho(R(z),x)+\rho(z,R(x))+\lambda\rho(z,x)))R_{1}(y)
\\&=&-R_{1}([x,y,z]_{R})+[R(x),R(y),R_{1}(z)]-R([R(x),y,R_{1}(z)]+[x,R(y),R_{1}(z)]\\&&+\lambda[x,y,R_{1}(z)])
+[R(y),R(z),R_{1}(x)]-R([R(y),z,R_{1}(x)]+[y,R(z),R_{1}(x)]\\&&+\lambda
[y,z,R_{1}(x)])
+[R(z),R(x),R_{1}(y)]-R([R(z),x,R_{1}(y)]+[z,R(x),R_{1}(y)]\\&&+\lambda
[z,x,R_{1}(y)]) ,
\end{eqnarray*}
and
\begin{eqnarray*}\delta \pi_{1}(x,y,z)&=&\pi_{1}(R(x),R(y),R(z))-\lambda^2 R \pi_{1}(x,y,z)-\lambda R
(\pi_{1}(Rx,y,z)+\pi_{1}(x,Ry,z)\\&&+\pi_{1}(x,y,Rz))-R(\pi_{1}(Rx,Ry,z)+\pi_{1}(x,Ry,Rz)+\pi_{1}(Rx,y,Rz))
 .
\end{eqnarray*}
Hence, $\partial_{R}R_1+\delta \pi_1=0$ is equivalent to (4.4) when
$n=1$.
\end{proof}

{\bf Definition 4.2.} A
 deformation of the Rota-Baxter 3-Lie algebra $(A,\pi,R
 )$ is said to be trivial if there is a linear map $N:A
 \longrightarrow A
 $ such that $K_{t}=I+tN~(\forall~ t)$ satisfies
 $$K_{t}R
 =R
 K_{t}$$ and
 $$K_{t}[x,y,z]_{t}=[K_{t}x, K_{t}y,K_{t}z]
 .$$

Based on the case of 3-Lie algebras, we give the
definition of a Nijenhuis operator of on Rota-Baxter 3-Lie algebras.\\

 {\bf Definition 4.3.} Let $(A,R)$
 be a Rota-Baxter 3-Lie algebra. A linear map $N:A
 \longrightarrow A $ is called a Nijenhuis operator if $NR
 =RN$ and $N$ is a Nijenhuis operator of $A$, that is,
\begin{eqnarray*}[Nx,Ny,Nz]
 &=&N([Nx,Ny,z]+[x,Ny,Nz]+[Nx,y,Nz])\\&&-N^{2}([Nx,y,z]+[x,Ny,z]+[x,y,Nz])+N^{3}([x,y,z])
 \end{eqnarray*}
for any $x,~y,~z \in A
 $.\\

 {\bf Remark 4.4.}
$(A,[ \ , \ , ]_{N},R
 )$ is also a Rota-Baxter 3-Lie algebra,
where
\begin{eqnarray*}[x,y,z]_{N}&=&[Nx,Ny,z]+[x,Ny,Nz]+[Nx,y,Nz]\\&&-N([Nx,y,z]+[x,Ny,z]+[x,y,Nz])+N^{2}([x,y,z]).
\end{eqnarray*}

{\bf Definition 4.5.} A linear map $K:V\longrightarrow A$
 is called an $\mathcal{O}$-operator on the Rota-Baxter 3-Lie algebra
 $A$ associated with the representation
$(V,\rho,R_{V})$ if $K:V\longrightarrow A$
 is an $\mathcal{O}$-operator on the 3-Lie algebra $A$ associated to the representation
$(V,\rho)$, that is, for all $u,v,w\in V$,
$$[K(u), K(v), K(w)]=K(\rho (K(u),K(v))w+\rho(K(v),K(w))u+\rho (K(w),K(u))v).$$
and $KR_{V
 }=RK$. Then $K$ is called an $\mathcal{O}$-operator on Rota-Baxter 3-Lie algebras.\\

{\bf Proposition 4.6.} Let $K:V\longrightarrow A$ be an
$\mathcal{O}$-operator on the Rota-Baxter 3-Lie algebra
 $(A,R)$ associated to the representation
$(V,\rho,R_{V})$. Then $(V,[ \ , \ , \ ]_{K},R_V)$ is Rota-Baxter
3-Lie algebra, where
 $$ [u,v,w]_{K}=\rho (K(u),K(v))w+\rho(K(v),K(w))u+\rho (K(w),K(u))v $$ for all $u,v,w \in
 V$.
\begin{proof} Due to Proposition 4.6 in \cite{35}, $(V,[ \ , \ , \ ]_{K})$ is a 3-Lie algebra.
In view of (2.3), by direct computation, we can check that $R_{V}$
is a Rota-Baxter operator of $(V,[ \, \, \ ]_{K})$.
\end{proof}

{\bf Proposition 4.7.} Let $K:V\longrightarrow A$ be an
$\mathcal{O}$-operator on Rota-Baxter 3-Lie algebra $(A,R)$
 associated to the representation $(V,\rho,R_{V})$. Define
 bilinear map
$\varrho_{K}:V\wedge V\longrightarrow \mathfrak{gl}(A)$ by
$$\varrho_{K}(u,v)x=[K(u),K(v),x]-K(\rho(Kv,x)u+\rho(x,Ku)v).$$
Then $(A,\varrho_{K},R)$ is a representation of $(V,[ \, \, \
]_{K},R_{V})$.
\begin{proof} In the light of Proposition 4.7 in \cite{35}, we know that $(A,\varrho_{K})$
is a representation of $(V,[ \, \, \ ]_{K})$. In the following, we
verify that (2.3) holds for $\varrho_{K}$. In fact,
\begin{eqnarray*}
&&\varrho_{K}(R_{V}u\wedge R_{V}v)R(x)\\
&=&[KR_{V}(u),KR_{V}(v),R(x)]-K\rho(KR_{V}v\wedge Rx
)R_{V}(u)-K\rho(Rx\wedge KR_{V}u)R_{V}v
\\
&=&R([RK(u),KR_{V}(v),x]+[RK(u),K(v),R(x)]+[K(u),RK(v),R(x)])\\&&+\lambda
R([RK(u),K(v),x]+[K(u),RK(v),x]+[K(u),K(v),R(x)]+\lambda
+[K(u),K(v),x]) \\&&-KR_{V}(\rho(RK(v)\wedge R(x)
)u+R_{V}(\rho(K(v), R(x)) +\rho( RK(v), x)+\lambda \rho(Kv,
x))R_{V}(u)\\&&+\lambda R_{V}(\rho(K(v), R(x))+\rho(RK(v),
x)+\lambda \rho(Kv,x ))u -KR_{V}(\rho(R(x),
RK(u))v\\&&-R_{V}(\rho(R(x), K(u))+\rho(x, RK(u))+\lambda \rho(x,
Ku))R_{V}(v)-\lambda R_{V}(\rho(R(x), K(u))\\&&+\rho(x,
RK(u))+\lambda\rho( x,Ku))v
\\&=&R\varrho_{K}(R_{V}(u),R_{V}(v))x+R(\varrho_{K}(R_{V}(u),
v)+\varrho_{K}(u, R_{V}(v))+\lambda \varrho_{K}(u,
v))R(x)\\&&+\lambda R(\varrho_{K}(R_{V}(u), v)+\varrho_{K}(u,
R_{V}(v))+\lambda \varrho_{K}(u, v))x,
\end{eqnarray*}
which follows that
 \begin{eqnarray*}\varrho_{K}(R_{V}u,
R_{V}v)R&=&R\varrho_{K}(R_{V}(u),R_{V}(v))+R(\varrho_{K}(R_{V}(u),
v)+\varrho_{K}(u, R_{V}(v))+\lambda \varrho_{K}(u, v))R\\&&+\lambda
R(\varrho_{K}(R_{V}(u), v)+\varrho_{K}(u, R_{V}(v))+\lambda
\varrho_{K}(u, v)).\end{eqnarray*}

 Hence, $(A,\varrho_{K},R)$ is a representation
of $(V,[ \ , \ , \ ]_{K},R_{V})$.
\end{proof}

It is natural to study the cohomology of the Rota-Baxter 3-Lie
algebra $(V,[ \ , \ , \ ]_{K},R_{V})$ with coefficients in the
representation $(A,\varrho_{K},R)$. We do not intend to discuss it
here and leave it to the readers.

\section{Extension of Rota-Baxter 3-Lie algebras}
Let $(A,R)$ be a Rota-Baxter 3-Lie algebra and $(V,R_{V})$ an
abelian Rota-Baxter 3-Lie algebra with the trivial product.\\

{\bf Definition 5.1.} A central extension of Rota-Baxter 3-Lie
algebra $(A,R)$ by the abelian Rota-Baxter 3-Lie algebra $(V,R_{V})$
is an exact sequence of Rota-Baxter 3-Lie algebras
$$0\longrightarrow(V,R_{V})\stackrel{i}{\longrightarrow} (\hat{A},R_{\hat{A}})\stackrel{p}{\longrightarrow} (A,R)\longrightarrow0$$
such that $[V,V,\hat{A}]=[V,\hat{A},V]=[\hat{A},V,V]=0$.\\

{\bf Definition 5.2.} Let $(\hat{A_1},R_{\hat{A_1}})$ and
$(\hat{A_2},R_{\hat{A_2}})$ be two central extensions of $(A,R)$ by
$(V,R_{V})$ . They are said to be equivalent if there is a
homomorphism of Rota-Baxter 3-Lie algebras
$\varphi:(\hat{A_1},R_{\hat{A_1}})\longrightarrow
(\hat{A_2},R_{\hat{A_2}})$ such that the following commutative
diagram holds:
$$\xymatrix{
  0 \ar[r] & (V,R_{V}) \ar[d]_{id} \ar[r]^{i} & (\hat{A_1},R_{\hat{A_1}} )\ar[d]_{\varphi} \ar[r]^{p} & (A,R) \ar[d]^{id} \ar[r] & 0\\
 0 \ar[r] & (V,R_{V}) \ar[r]^{i} & (\hat{A_2},R_{\hat{A_2}}) \ar[r]^{p} & (A,R)  \ar[r] & 0
 .}$$

Similarly to the case of 3-Lie algebra \cite{37,38,39}, we easily
get the
following result:\\

{\bf Lemma 5.3.} With the above notions, $(V,\rho, R_{V})$ is a
representation of $(A,R)$ and is independent on the choice of the
section $s$. Moreover, equivalent abelian extensions give the same
representation.\\

A section of a central extension $(\hat{A},R_{\hat{A}})$ of $(A,R)$
by $(V,R_{V})$ is a linear map $s:A\longrightarrow \hat{A}$ such
that $ps=I$.

Let $s:A\longrightarrow \hat{A}$ be any section of $p$. Define two
multilinear map $\psi:\wedge^{3}A\longrightarrow V$ and linear map
$\chi:A\longrightarrow V$ by
$$\psi(x,y,z)=[s(x),s(y),s(z)]-s[x,y,z]$$ and
$$\chi(x)=R_{\hat{A}}(s(x))-s(R(x)).$$

{\bf Proposition 5.4.} The vector space $(A\oplus V,[ \ , \ , \
]_{\psi},R_{\chi})$ with the bracket
$$[x+a,y+b,z+c]_{\psi}=[x,y,z]+\psi(x,y,z)$$
and
$$R_{\chi}(x+a)=R(x)+R_{V}(a)+\chi(x)$$
is a Rota-Baxter 3-Lie algebra if and only if $(\psi,\chi)$ is a
2-cocycle of $(A,R)$ with coefficients in the trivial representation
$(V,R_{V})$.

\begin{proof}
If $(\psi,\chi)\in
C_{\hbox{RB}}^{2}(A,V)=C_{\hbox{3-Lie}}^{2}(A,V)\times
C_{\hbox{3-Lie}}^{1}(A_R,\tilde{V}) $ is a 2-cocycle, then
$$\partial_{RB}(\psi,\chi)=(\partial\psi,\partial_{R}\chi+(-1)^{2}\delta\psi)=0.$$

Based on deformation theory of 3-Lie algebras \cite{37,38,39}, we
only need to check that $R_{\chi}$ is a Rota-Baxter operator of the
3-Lie algebra $A\oplus V$ if and only if
$\partial_{R}\chi+\delta\psi=0$.

On the one hand, by computation,
\begin{eqnarray*}&&R_{\chi}([R_{\chi}(x+a),R_{\chi}(y+b),z+c]_{\psi}+[R_{\chi}(x+a),y+b,R_{\chi}(z+c)]_{\psi} \\&&+[x+a,R_{\chi}(y+b),R_{\chi}(z+c)]_{\psi}
 +\lambda([R_{\chi}(x+a),y+b,z+c]_{\psi}+[x+a,R_{\chi}(y+b),z+c]_{\psi} \\&&+[x+a,y+b,R_{\chi}(z+c)]_{\psi}
 +\lambda[x+a,y+b,z+c]_{\psi}))-[R_{\chi}(x+a),R_{\chi}(y+b),R_{\chi}(z+c)]_{\psi}\\&=&
R([Rx,Ry,z]+[Rx,y,Rz]+[x,Ry,Rz]+\lambda([Rx,y,z]+[x,Ry,z]+[x,y,Rz])\\&&+\lambda^{2}[x,y,z])
+R_{V}(\psi(Rx,Ry,z)+\psi(Rx,y,Rz)+\psi(x,Ry,Rz)
+\lambda(\psi(Rx,y,z)\\&&+\psi(x,Ry,z)+\psi(x,y,Rz))
+\lambda^{2}\psi(x,y,z))-[Rx,Ry,Rz]-\psi(Rx,Ry,Rz).~~~~~~(5.1)
 \end{eqnarray*}
On the other hand, since $(V,R_{V})$ is a trivial representation,
\begin{eqnarray*}&&
(\partial_{R}\chi+\delta\psi)(x,y,z)
\\&=&-\chi([x,y,z]_{R})+\psi(Rx,Ry,Rz)-R_{V}\lambda^{2}\psi(x,y,z)-R_{V}\lambda(\psi(Rx,y,z)+\psi(x,Ry,z)\\&&+\psi(x,y,Rz))
-R_{V}(\psi(Rx,Ry,Rz)+\psi(Rx,Ry,z)+\psi(Rx,y,Rz)).~~~~~~~~~~~~~~~~~(5.2)
 \end{eqnarray*}
In the light of (5.1) and (5.2), we get the conclusion.
\end{proof}
\section{Appendix: Proof of Proposition 3.2}
\begin{proof} Denote $f_{i}^{k}$ by
$$f_{i}^{k}(x_1,\cdot\cdot\cdot,x_n)=
f(\underbrace{I\otimes\cdot\cdot\cdot \otimes
R\otimes\cdot\cdot\cdot \otimes I}_{R~~ \hbox{appears} ~~k
~~\hbox{times  but not in the ith
place}})(x_1,\cdot\cdot\cdot,x_n),$$ and $\tilde{R}=R\wedge R$.

For any $f\in C_{\hbox{3-Lie}}^{n}(A,V)$ and $X_i=x_{i}\wedge
y_{i}\in L $ and $z\in A$. For convenient, we denote
$(X_1,\cdot\cdot\cdot,X_{n},z)=(x_1,x_2,\cdot\cdot\cdot,x_{2n+1})\in
L^{\otimes n}\wedge A$.

On the one hand, by direct calculation,
\begin{eqnarray*}\delta \partial f (X_1,\cdot\cdot\cdot,X_{n},z)&=&
\partial
f(\tilde{R}(X_1),\cdot\cdot\cdot,\tilde{R}(X_{n}),R(z))-R_{V}\sum_{k=0}^{2n}\lambda^{2n-k}(\partial
f)^{k}(X_1,\cdot\cdot\cdot,X_{n},z)\\&=&B_1+B_2+B_3+B_4+B_5-C_1-C_2-C_3-C_4-C_5,
\end{eqnarray*}
with
\begin{eqnarray*}B_1&=&\sum_{1\leq i< k\leq
{n}}
(-1)^{i}f(R(x_1),\cdot\cdot\cdot,R(x_{2j-1}),[R(x_{2i-1}),R(x_{2i}),R(x_{2j})],
R(x_{2j+1}),\cdot\cdot\cdot,R(x_{2n+1}))
\\&&+\sum_{1\leq i< j\leq
{n}}
(-1)^{i}f(R(x_1),\cdot\cdot\cdot,R(x_{2j-2}),[R(x_{2i-1}),R(x_{2i}),R(x_{2j-1})],
R(x_{2j}),\cdot\cdot\cdot,R(x_{2n+1})),\end{eqnarray*}
$$B_2=
\sum_{i=1}^{n}(-1)^{i}f(R(x_1),\cdot\cdot\cdot,R(x_{2n}),[R(x_{2i-1}),R(x_{2i}),R(x_{2n+1})]),$$
$$B_3=
\sum_{i=1}^{n}(-1)^{i+1}\rho(R(x_{2i-1}),R(x_{2i}))f(R(x_1),\cdot\cdot\cdot,R(x_{2n+1})),$$
\begin{eqnarray*}B_4&=&
(-1)^{n+1}(\rho(R(x_{2n}),R(x_{2n+1}))f(R(x_1),\cdot\cdot\cdot,R(x_{2n-2}),R(x_{2n-1}))
\\&&+\rho(R(x_{2n+1}),R(x_{2n-1})f(R(x_1),\cdot\cdot\cdot,R(x_{2n-2}),R(x_{2n}))),
\end{eqnarray*}
\begin{eqnarray*}C_1
&=&R_{V}\sum_{1\leq i< j\leq {n}}
(-1)^{i}\sum_{k=0}^{2n-2}(\lambda^{2n-k}f_{2j-2}^{k}(x_1,\cdot\cdot\cdot,x_{2j-1},[x_{2i-1},x_{2i},x_{2j}],
\cdot\cdot\cdot,x_{2n+1})
\\&&+\sum_{k=0}^{2n-2}\lambda^{2n-k-1}f_{2j-2}^{k}(x_1,\cdot\cdot\cdot,x_{2j-1},[R(x_{2i-1}),x_{2i},x_{2j}],
\cdot\cdot\cdot,x_{2n+1})
\\&&+\sum_{k=0}^{2n-2}\lambda^{2n-k-1}f_{2j-2}^{k}(x_1,\cdot\cdot\cdot,x_{2j-1},[x_{2i-1},R(x_{2i}),x_{2j}],
\cdot\cdot\cdot,x_{2n+1})
\\&&+\sum_{k=0}^{2n-2}\lambda^{2n-k-1}f_{2j-2}^{k}(x_1,\cdot\cdot\cdot,x_{2j-1},[x_{2i-1},x_{2i},R(x_{2j})],
\cdot\cdot\cdot,x_{2n+1})
\\&&+\sum_{k=0}^{2n-2}\lambda^{2n-k-2}f_{2j-2}^{k}(x_1,\cdot\cdot\cdot,x_{2j-1},[R(x_{2i-1}),R(x_{2i}),x_{2j}],
\cdot\cdot\cdot,x_{2n+1})
\\&&+\sum_{k=0}^{2n-2}\lambda^{2n-k-2}f_{2j-2}^{k}(x_1,\cdot\cdot\cdot,x_{2j-1},[R(x_{2i-1}),x_{2i},R(x_{2j})],
\cdot\cdot\cdot,x_{2n+1})
\\&&+\sum_{k=0}^{2n-2}\lambda^{2n-k-2}f_{2j-2}^{k}(x_1,\cdot\cdot\cdot,x_{2j-1},[x_{2i-1},R(x_{2i}),R(x_{2j})],
\cdot\cdot\cdot,x_{2n+1})
\\&&+\sum_{k=0}^{2n-3}\lambda^{2n-k-3}f_{2j-2}^{k}(x_1,\cdot\cdot\cdot,x_{2j-1},[R(x_{2i-1}),R(x_{2i}),R(x_{2j})],
\cdot\cdot\cdot,x_{2n+1})),
\end{eqnarray*}
\begin{eqnarray*}C_2
&=&R_{V}\sum_{1\leq i< j\leq {n}}
(-1)^{i}(\sum_{k=0}^{2n-2}\lambda^{2n-k}f_{2j-3}^{k}(x_1,\cdot\cdot\cdot,x_{2j-2},[x_{2i-1},x_{2i},x_{2j-1}],x_{2j}
\cdot\cdot\cdot,x_{2n+1})
\\&&+\sum_{k=0}^{2n-2}\lambda^{2n-k-1}f_{2j-3}^{k}(x_1,\cdot\cdot\cdot,x_{2j-2},[R(x_{2i-1}),x_{2i},x_{2j-1}],x_{2j}
\cdot\cdot\cdot,x_{2n+1})
\\&&+\sum_{k=0}^{2n-2}\lambda^{2n-k-1}f_{2j-3}^{k}(x_1,\cdot\cdot\cdot,x_{2j-2},[x_{2i-1},R(x_{2i}),x_{2j-1}],x_{2j}
\cdot\cdot\cdot,x_{2n+1})
\\&&+\sum_{k=0}^{2n-2}\lambda^{2n-k-1}f_{2j-3}^{k}(x_1,\cdot\cdot\cdot,x_{2j-2},[x_{2i-1},x_{2i},R(x_{2j-1})],x_{2j}
\cdot\cdot\cdot,x_{2n+1})
\\&&+\sum_{k=0}^{2n-2}\lambda^{2n-k-2}f_{2j-3}^{k}(x_1,\cdot\cdot\cdot,x_{2j-2},[R(x_{2i-1}),R(x_{2i}),x_{2j-1}],x_{2j}
\cdot\cdot\cdot,x_{2n+1})
\\&&+\sum_{k=0}^{2n-2}\lambda^{2n-k-2}f_{2j-3}^{k}(x_1,\cdot\cdot\cdot,x_{2j-2},[R(x_{2i-1}),x_{2i},R(x_{2j-1})],x_{2j}
\cdot\cdot\cdot,x_{2n+1})
\\&&+\sum_{k=0}^{2n-2}\lambda^{2n-k-2}f_{2j-3}^{k}(x_1,\cdot\cdot\cdot,x_{2j-2},[x_{2i-1},R(x_{2i}),R(x_{2j-1})],x_{2j}
\cdot\cdot\cdot,x_{2n+1})
\\&&+\sum_{k=0}^{2n-3}\lambda^{2n-k-3}f_{2j-3}^{k}(x_1,\cdot\cdot\cdot,x_{2j-2},[R(x_{2i-1}),R(x_{2i}),R(x_{2j-1})],x_{2j}
\cdot\cdot\cdot,x_{2n+1})),
\end{eqnarray*}
\begin{eqnarray*}C_3
&=&R_{V}\sum_{i=1}^{n}(-1)^{i}(\sum_{k=0}^{2n-2}\lambda^{2n-k}f_{2n-1}^{k}(x_1,\cdot\cdot\cdot,x_{2n},[x_{2i-1},x_{2i},x_{2n+1}])
\\&&+\sum_{k=0}^{2n-2}\lambda^{2n-1-k}f_{2n-1}^{k}(x_1,\cdot\cdot\cdot,x_{2n},[R(x_{2i-1}),x_{2i},z])
\\&&+\sum_{k=0}^{2n-2}\lambda^{2n-1-k}f_{2n-1}^{k}(x_1,\cdot\cdot\cdot,x_{2n},[x_{2i-1},R(x_{2i}),z])
\\&&+\sum_{k=0}^{2n-2}\lambda^{2n-1-k}f_{2n-1}^{k}(x_1,\cdot\cdot\cdot,x_{2n},[x_{2i-1},x_{2i},R(z)])
\\&&+\sum_{k=0}^{2n-2}\lambda^{2n-2-k}f_{2n-1}^{k}(x_1,\cdot\cdot\cdot,x_{2n},[R(x_{2i-1}),R(x_{2i}),z])
\\&&+\sum_{k=0}^{2n-2}\lambda^{2n-2-k}f_{2n-1}^{k}(x_1,\cdot\cdot\cdot,x_{2n},[R(x_{2i-1}),x_{2i},R(z)])
\\&&+\sum_{k=0}^{2n-2}\lambda^{2n-2-k}f_{2n-1}^{k}(x_1,\cdot\cdot\cdot,x_{2n},[x_{2i-1},R(x_{2i}),R(z)])
\\&&+\sum_{k=0}^{2n-3}\lambda^{2n-3-k}f_{2n-1}^{k}(x_1,\cdot\cdot\cdot,x_{2n},[R(x_{2i-1}),R(x_{2i}),R(z)]),
\end{eqnarray*}
\begin{eqnarray*}C_4
&=&R_{V}\sum_{i=1}^{n}(-1)^{i+1}(\sum_{k=0}^{2n-1}\lambda^{2n-k}\rho(x_{2i-1},x_{2i})f^{k}(x_1,\cdot\cdot\cdot,x_{2n},x_{2n+1})
\\&&
+\sum_{k=0}^{2n-1}\lambda^{2n-1-k}\rho(R(x_{2i-1}),x_{2i})f^{k}(x_1,\cdot\cdot\cdot,x_{2n},x_{2n+1})
\\&&
+\sum_{k=0}^{2n-1}\lambda^{2n-1-k}\rho(x_{2i-1},R(x_{2i}))f^{k}(x_1,\cdot\cdot\cdot,x_{2n},x_{2n+1})
\\&&
+\sum_{k=0}^{2n-2}\lambda^{2n-2-k}\rho(R(x_{2i-1}),R(x_{2i}))f^{k}(x_1,\cdot\cdot\cdot,x_{2n},x_{2n+1})),\end{eqnarray*}
and
\begin{eqnarray*}C_5
&=&R_{V}
(-1)^{n+1}(\sum_{k=0}^{2n-1}\lambda^{2n-k}\rho(x_{2n},x_{2n+1})f^{k}(x_1,\cdot\cdot\cdot,x_{2n-2},x_{2n-1})
\\&&+\sum_{k=0}^{2n-1}\lambda^{2n-1-k}\rho(R(x_{2n}),x_{2n+1})f^{k}(x_1,\cdot\cdot\cdot,x_{2n-2},x_{2n-1})
\\&&+\sum_{k=0}^{2n-1}\lambda^{2n-1-k}\rho(x_{2n},R(x_{2n+1}))f^{k}(x_1,\cdot\cdot\cdot,x_{2n-2},x_{2n-1})
\\&&+\sum_{k=0}^{2n-2}\lambda^{2n-2-k}\rho(R(x_{2n}),R(x_{2n+1}))f^{k}(x_1,\cdot\cdot\cdot,x_{2n-2},x_{2n-1})
\\&&+\sum_{k=0}^{2n-1}\lambda^{2n-k}\rho(x_{2n+1},x_{2n-1})f(x_1,\cdot\cdot\cdot,x_{2n-2},x_{2n})
\\&&+\sum_{k=0}^{2n-1}\lambda^{2n-1-k}\rho(R(x_{2n+1}),x_{2n-1})f(x_1,\cdot\cdot\cdot,x_{2n-2},x_{2n})
\\&&+\sum_{k=0}^{2n-1}\lambda^{2n-1-k}\rho(x_{2n+1},R(x_{2n-1}))f(x_1,\cdot\cdot\cdot,x_{2n-2},x_{2n})
\\&&+\sum_{k=0}^{2n-2}\lambda^{2n-2-k}\rho(R(x_{2n+1}),R(x_{2n-1}))f(x_1,\cdot\cdot\cdot,x_{2n-2},x_{2n}))
.
\end{eqnarray*}
On the other hand, denote
$$\partial_{R}\delta
f(X_1,\cdot\cdot\cdot,X_{n},z)=D_1-D_2+E+F+G,$$ where
\begin{eqnarray*}D_1&=&
\sum_{1\leq i< j\leq {n}} (-1)^{i}
f(\tilde{R}(X_1),\cdot\cdot\cdot,\tilde{R}(X_{j-1}),\tilde{R}([X_i,X_j]_F),\tilde{R}(X_{j+1}),\cdot\cdot\cdot,\tilde{R}(X_{n})
,R(z))\\ &=&\sum_{1\leq i< j\leq {n}} (-1)^{i}
f(R(x_1),\cdot\cdot\cdot,R(x_{2j-2}),R(x_{2j-1}),R[x_{2i-1},x_{2i},x_{2j}]_R,R(x_{2j+1}),\cdot\cdot\cdot,R(x_{2n+1})
)
\\&&+\sum_{1\leq i< j\leq {n}} (-1)^{i}
f(R(x_1),\cdot\cdot\cdot,R(x_{2j-2}),R[x_{2i-1},x_{2i},x_{2j-1}]_R,R(x_{2j}),\cdot\cdot\cdot,R(x_{2n+1})
),
\end{eqnarray*}
\begin{eqnarray*}D_2&=&R_{V}\sum_{k=0}^{2n-2}\lambda^{2n-2-k}f^{k}(X_1,\cdot\cdot\cdot,\hat{X_i},X_{j-1},
[X_i,X_j]_F,X_{j+1},\cdot\cdot\cdot,X_{n}
,z)\\
&=&R_{V}\sum_{k=0}^{2n-2}\lambda^{2n-2-k}f^{k}(x_1,\cdot\cdot\cdot,x_{2j-3},x_{2j-2},x_{2j-1},
[x_{2i-1},x_{2i},x_{2j}]_R,x_{2j+1},\cdot\cdot\cdot,x_{2n+1}
)\\&&+R_{V}\sum_{k=0}^{2n-2}\lambda^{2n-2-k}f^{k}(x_1,\cdot\cdot\cdot,x_{2j-3},x_{2j-2},
[x_{2i-1},x_{2i},x_{2j-1}]_R,x_{2j},x_{2j+1},\cdot\cdot\cdot,x_{2n+1}
)
\\
&=&R_{V}\sum_{k=0}^{2n-2}\lambda^{2n-2-k}f^{k}_{2j-2}(x_1,\cdot\cdot\cdot,x_{2j-3},x_{2j-2},x_{2j-1},
[x_{2i-1},x_{2i},x_{2j}]_R,x_{2j+1},\cdot\cdot\cdot,x_{2n+1}
)\\&&+R_{V}\sum_{k=0}^{2n-3}\lambda^{2n-3-k}f^{k}_{2j-2}(x_1,\cdot\cdot\cdot,x_{2j-3},x_{2j-2},x_{2j-1},
R[x_{2i-1},x_{2i},x_{2j}]_R,x_{2j+1},\cdot\cdot\cdot,x_{2n+1} )
\\&&+R_{V}\sum_{k=0}^{2n-2}\lambda^{2n-2-k}f^{k}_{2j-3}(x_1,\cdot\cdot\cdot,x_{2j-3},x_{2j-2},
[x_{2i-1},x_{2i},x_{2j-1}]_R,x_{2j},x_{2j+1},\cdot\cdot\cdot,x_{2n+1}
)
\\&&+R_{V}\sum_{k=0}^{2n-3}\lambda^{2n-3-k}f^{k}_{2j-3}(x_1,\cdot\cdot\cdot,x_{2j-3},x_{2j-2},
R[x_{2i-1},x_{2i},x_{2j-1}]_R,x_{2j},x_{2j+1},\cdot\cdot\cdot,x_{2n+1}
),\end{eqnarray*}
\begin{eqnarray*}E&=&\sum_{i=1}^{n}(-1)^{i}\delta
f(X_1,\cdot\cdot\cdot,\hat{X_i},\cdot\cdot\cdot,X_{n},[X_i,z]_{R})\\&=&\sum_{i=1}^{n}(-1)^{i}f(R(x_1),\cdot\cdot\cdot,
R(x_{2n}),R[x_{2i-1},x_{2i},x_{2n+1}]_R
)\\&&-R_{V}\sum_{k=0}^{2n-2}\lambda^{2n-2-k}f^{k}(x_1,\cdot\cdot\cdot,
x_{2n},[x_{2i-1},x_{2i},x_{2n+1}]_R )
\\&=&\sum_{i=1}^{n}(-1)^{i}f(R(x_1),\cdot\cdot\cdot,
R(x_{2n}),R[x_{2i-1},x_{2i},x_{2n+1}]_R
)\\&&-R_{V}\sum_{k=0}^{2n-2}\lambda^{2n-2-k}f_{2n-1}^{k}(x_1,\cdot\cdot\cdot,
x_{2n},[x_{2i-1},x_{2i},x_{2n+1}]_R
)\\&&-R_{V}\sum_{k=0}^{2n-3}\lambda^{2n-3-k}f_{2n-1}^{k}(x_1,\cdot\cdot\cdot,
x_{2n},R[x_{2i-1},x_{2i},x_{2n+1}]_R )
\end{eqnarray*}
\begin{eqnarray*}F&=&
\sum_{i=1}^{n}(-1)^{i+1}\tilde{\rho}(X_i)\delta
f(X_1,\cdot\cdot\cdot,\hat{X_i},\cdot\cdot\cdot,X_{n},z)\\&=&\sum_{i=1}^{n}(-1)^{i+1}\tilde{\rho}(X_i)(
f(R(x_1),\cdot\cdot\cdot, R(x_{2n}),R(x_{2n+1}) )
-R_{V}\sum_{k=0}^{2n-2}\lambda^{2n-2-k}f^{k}(x_1,\cdot\cdot\cdot,
x_{2n+1} )),
\end{eqnarray*}
and
\begin{eqnarray*}G&=&
(-1)^{n+1}(\tilde{\rho}(x_{2n},x_{2n+1})\delta
f(x_1,\cdot\cdot\cdot,x_{2n-2},x_{2n-1})+\tilde{\rho}(x_{2n+1},x_{2n-1})\delta
f(x_1,\cdot\cdot\cdot,x_{2n-2},x_{2n}))
\\&=&(-1)^{n+1}\tilde{\rho}(x_{2n},x_{2n+1})(f(R(x_1),\cdot\cdot\cdot,
R(x_{2n-2}),R(x_{2n-1})-R_{V}\sum_{k=0}^{2n-2}\lambda^{2n-2-k}\times
\\&&f^{k}(x_1,\cdot\cdot\cdot,
x_{2n-2},x_{2n-1})+(-1)^{n+1}\tilde{\rho}(x_{2n+1},x_{2n-1})(f(R(x_1),\cdot\cdot\cdot,
R(x_{2n-2}),R(x_{2n})\\&&-R_{V}\sum_{k=0}^{2n-2}\lambda^{2n-2-k}f^{k}(x_1,\cdot\cdot\cdot,
x_{2n-2},x_{2n}).
\end{eqnarray*}
Using 2.1 and 2.2, obviously, $B_1-C_1-C_2=D_1-D_2$. Similarly,
$B_2-C_3=E$, $B_3-C_4=F$ and $B_4-C_5=G$. These yield that
$\partial_{R}\delta=\delta
\partial$.
\end{proof}

\begin{center}{\textbf{Acknowledgments}}
\end{center}
Project supported by the National Natural Science Foundation of
China ( No. 11871421, No. 11401530 ) and the Natural Science
Foundation of Zhejiang Province of China ( No. LY19A010001 ).



\begin{thebibliography} {99}

\bibitem [1] {1} M. Aguiar, Pre-Poisson algebras, Lett. Math. Phys. 54 (2000),
263-277.
\bibitem [2] {2} J. A. de Azc'arraga, J. M. Izquierdo, Cohomology of Filippov algebras
and an analogue of Whitehead's lemma, J. Phys.: Conf. Ser. 175
(2019), 012001.
\bibitem [3] {3} C. Bai, A unified algebraic approach to classical Yang-Baxter
equation, J. Phys. A Math. Theor. 40(2007), 11073-11082.
\bibitem [4] {4} D. Balavoine, Deformation of algebras over a quadratic operad, Contemp. Math. AMS, 202 (1997), 207-234.
\bibitem [5]{5} G. Baxter, An analytic problem whose solution follows from a simple
algebraic identity, Pac. J. Math. 10 (1960), 731-742.
\bibitem [6]{6} C. Bai, O. Bellier, L. Guo, X. Ni, Spliting of operations, Manin
products and Rota-Baxter operators, Int. Math. Res. Not. 2013(3),
485-524.
\bibitem [7]{7} C. M. Bai, L. Guo, Y. Sheng, Bialgebras, the classical Yang-Baxter equation and Manin triples for 3-Lie algebras,
 arXiv:1604.05996.
\bibitem [8]{8} R. P. Bai, L. Guo, J. Q. Li, Y. Wu, Rota-Baxter 3-Lie algebras, J. Math. Phys. 54
(2013), 063504.
\bibitem [9]{9}A. Connes, D. Kreimer, Renormalization in quantum field theory
and the Riemann-Hilbert problem. I. The Hopf algebra structure of
graphs and the main theorem, Comm. Math. Phys. 210 (2000) 249-273.
\bibitem[10]{10} Y. Daletskii, L. Takhtajan, Leibniz and Lie algebra structures
for Nambu algebra, Lett. Math. Phys., 39 (1997), 127-141.
\bibitem [11]{11} A. Das, Deformations of associative Rota-Baxter operators, J.
Algebra 560 (2020), 144-180.
\bibitem [12]{12} A. Das, Cohomology and deformations of twisted Rota-Baxter operators
and NS-algebras, arXiv: 2010.01156.
\bibitem [13]{13} A. Das, S. Guo, Twisted relative Rota-Baxter operators on Leibniz
algebras and NS-Leibniz algebras, arXiv: 2102.09752.
\bibitem [14]{14} A. Das, Cohomology and deformations of weighted Rota-Baxter
operators, arXiv: 2108.02627.
\bibitem [15]{15} A. Das, Cohomology of weighted Rota-Baxter Lie algebras and Rota-Baxter
paired operators, arXiv: 2109.01972v1.
\bibitem [16]{16} A. Das, S. K. Misha, The $L_{\infty}$-deformations of associative Rota-Baxter
algebras and homotopy Rota-Baxter operators, arXiv: 2008.11076.
\bibitem[17]{80} V. T. Filippov, n-Lie algebras, Sibirsk. Mat. Zh. 26 (1985), 126-140.
\bibitem[18]{17} P. Gautheron, Some remarks concerning Nambu mechanics, Lett. Math.
Phys., 37 (1996), 103-116.
\bibitem [19]{18} M. Gerstenhaber, The cohomology structure of an associative
 ring, Ann. Math. 78 (1963), 267-288.
\bibitem [20]{19} M. Gerstenhaber, On the deformation
of rings and algebras, Ann. Math. 79 (2) (1964), 59-103.
\bibitem [21]{20} N. S. Gu, L. Guo, Generating functions from the viewpoint of
Rota-Baxter algebras, Discrete Math. 338 (2015), 536-554.
\bibitem [22]{21} L. Guo, H. Lang, Y. Sheng, Integration and geometrization of
Rota-Baxter Lie algebras, Adv. Math. 387 (2021).
\bibitem [23]{22} L. Guo, W. Keigher, Baxter algebras and shuffle products, Adv. Math. 150 (2000), 117-149.
\bibitem [24]{23} L. Guo, Y. N. Li, Y. H. Sheng, G. D. Zhou, Cohomologies,
extensions and deformations of differential algebras of any weights,
arXiv:2003.03899v1 (2020).
\bibitem [25]{24} L. Guo, B. Zhang, Renormalization of multiple zeta values, J. Algebra 319 (2008), 3770-3809.
\bibitem [26]{25} J. Jiang, Y. Sheng, C. Zhu, Cohomologies of relative Rota-Baxter operators on Lie groups and Lie
algebras, arXiv: 2108.02627.
\bibitem [27]{26} H. L. Lang, Y. Sheng, Factorizable Lie bialgebras, quadratic Rota-Baxter Lie algebras and
Rota-Baxter Lie bialgebras, arXiv:2112.07902.
\bibitem [28]{27} A. Lazarev, Y. Sheng, R. Tang, Deformations and homotopy Theory
of Relative Rota-Baxter Lie Algebras, Comm. Math. Phys. 383 (1)
(2021), 595-631.
\bibitem [29]{28} A. Nijenhuis, R. Richardson, Cohomology and deformations in
graded Lie algebras, Bull. Amer. Math. Soc. 72 (2) (1966), 1-29.
 \bibitem [30]{29} A. Nijenhuis, R. Richardson, Commutative algebra cohomology and
deformations of Lie and associative algebras, J. Algebra 9 (1968),
42-105.
\bibitem [31] {30} M. Rotkiewicz, Cohomology ring of $n$-Lie algebras, Extracta Math.
20 (2005), 219-232.
\bibitem [32]{31} Q. X. Sun, Representations and cohomolgies of differential 3-Lie algebras with
any weights, preprint.
 \bibitem [33]{32} Q. X. Sun, S. Chen, Representations and Cohomologies
 of differential Lie-Yamaguti algebras with any weights, preprint.
\bibitem[34]{33} L. Takhtajan, A higher order analog of Chevally-Eilenberg complex and deformation
theory of n-algebras, St. Petersburg Math. J., 6 (1995), 429-438.
\bibitem [35]{34} R. Tang, C. Bai, L. Guo, Y. Sheng, Deformations and their
controlling cohomologies of $\mathcal{O}$-operators, Comm. Math.
Phys. 368 (2) (2019), 665-700.
\bibitem [36] {35} R. Tang, S. Hou, Y. H. Sheng, Lie 3-algebras
and deformationsof relative Rota-Baxter operators on 3-Lie algebras,
J. Alg. 567 (2021), 37-62.
\bibitem [37]{36}  K. Wang, G. D. Zhou, Deformations and homotopy theory of Rota-Baxter algebras of any
weight, arXiv: 2108.06744.
\bibitem [38] {37} R. R. Xu, Cohomology, derivations and abelian extensions of 3-Lie
algebras, J Algebra Appl, 2019.
\bibitem [39] {38} T. Zhang, Deformations and extensions of 3-Lie algebras, arXiv:1401.4656v4 (2020).
\bibitem [40] {39} T. Zhang, Cohomology and deformations of 3-Lie colour algebras,
Linear Multilinear A 63(2015), 651-671.




\end{thebibliography}
\end{document}